\documentclass{amsart}

\usepackage{amssymb,latexsym,amsxtra,amscd}
\usepackage[mathscr]{eucal}
\usepackage{enumerate}

\numberwithin{equation}{section}

\newtheorem{thm}{Theorem}[section]

\theoremstyle{definition}

\theoremstyle{remark}

\newtheorem{claim}{Claim}

\begin{document}

\title{On the (non)existence of states on orthogonally closed subspaces in
an inner product space}

\author{E.Chetcuti}
\address{
Emanuel Chetcuti \\
Department of Mathematics \\
University of Malta \\
Msida MSD.06, Malta}
\email {emanria@maltanet.net}
\author{P.Pt\'{a}k}
\thanks {The second author acknowledges the support of the Center of Machine
Perception, Czech Technical University and the grant no.\
201/02/1540 of the Grant Agency of the Czech Republic. The support enabled
the second author to work with the first author on the problem studied in
this paper.}
\address{
Pavel Pt\'{a}k \\
Department of Mathematics \\
Faculty of Electrical Engineering \\
Czech Technical University\\
166 27 Prague 6, Czech Republic}
\email {ptak@math.feld.cvut.cz}

\date{\today}



\begin{abstract}
Suppose that $S$ is an \emph {incomplete} inner product space. In~\cite
{Dvu} A.~Dvure\v{c}enskij shows that there are no finitely additive states
on orthogonally closed subspaces, $F(S)$, of $S$ that are regular with
respect to finitely dimensional spaces. In this note we show that the most
important special case of the former result---the case of the evaluations
given by vectors in the ``Gleason manner''---allows for a relatively simple
proof. This result further reinforces the conjecture that there are no
finitely additive states on $F(S)$ at all.
\end{abstract}

\maketitle

\section{Introduction}

Let $S$ be a real or complex separable inner product space and let
$\langle\cdot,\cdot\rangle$ denote the inner product of $S$. Let us denote
by $F(S)$ the set of all orthogonally closed subspaces of $S$. A subspace
$M$ of $S$ is in $F(S)$ if $M=M^{\perp\perp}$, where $M^\perp = \{x\in
S:\langle x,y\rangle=0 \text{ for all } y\in M\}$. It turns out that if we
understand $F(S)$ with the ordering given by the inclusion relation and with
orthocomplementation relation $M \longrightarrow M^\perp$ as defined above,
then $F(S)$ becomes a complete orthocomplemented lattice. However, $F(S)$
does not have to be orthomodular. In fact, Amemiya and Araki
\cite{Amem-Arak} proved the following algebraic criterion for the
(topological) completeness of an inner product space $S$: an inner product
space $S$ is complete if and only if $F(S)$ is orthomodular.

Let us now turn to measure-theoretic criteria for the completeness of $S$.
The following result by Hamhalter and Pt\'{a}k \cite{Ham-Ptak} initiated a
series of interesting measure theoretic characterizations for the
completeness of an inner product space $S$~\cite {Dvu}.

\begin{thm}
An inner product space $S$ is complete if and only if $F(S)$
possesses a $\sigma$-additive state.
\end{thm}

In 1988 \cite{Ptak}, Pt\'{a}k asked whether $S$ has to be complete if $F(S)$
possesses a finitely additive state. Recently, Dvure\v{c}enskij  and
Pt\'{a}k \cite{Dvu-Ptak} proved that if $S$ is an incomplete inner product
space, then the assumption that there is a finitely additive state on $F(S)$
implies that the range of this state has to be the entire interval $[0,1]$.
In this note we show that an inner product space $S$ is complete if, and
only if, there exists $u\in \overline{S}$ such that $s_u$ defines a state on
$F(S)$, where by $\overline{S}$ is denoted the completion of $S$. Here, for
any vector $u\in S$ with $\|u\|=1$, by $s_u$ is meant the ``Gleason''
assignment defined by
\begin{align*}
s_u:F(S)&\to[0,1]\\
M&\mapsto\langle P_{\overline M}u,u\rangle.
\end{align*}

Before we launch on the proof proper, let us summarize the ``state of the
art'' of the state problem for $F(S)$. If there are states on $F(S)$ then
there are pure states on $F(S)$ (Krein--Milman). But in view of the previous
two facts these pure states must be rather bizarre. Thus, a conjecture
remains that for an incomplete space $S$ the lattice $F(S)$ is stateless.

\section{Results}

Let $S$ be a separable inner product space and let $\overline{S}$ be its
completion. In this section we mainly prove the result formulated in the
introduction.

\begin{thm}\label{t:1} A separable inner product space $S$ is complete if,
and only if, there exists $u\in \overline{S}$ such that
 $$s_u: M\mapsto\langle P_{\overline M}u,u\rangle$$
defines a state on $F(S)$.
\end{thm}
\begin{proof} If $S$ is complete then, obviously, for every $u\in S =
\overline{S}$, $s_u$ is a ($\sigma-$additive) state on $F(S)$ ($F(S)=L(S)$
and this follows from Gleason's theorem).

For the second implication, suppose that there exists a vector $u\in
\overline{S}$ such that $s_u$ is a state on $F(S)$. We divide the proof
into auxiliary results. We believe that they could be of certain importance
in their own right.

\begin{claim}Suppose that there exists $u\in S$ such that
 $$s_u: M\mapsto\langle P_{\overline M}u,u\rangle$$
defines a state on $F(s)$. Then for every unit vector $v\in S$,
$s_v$ defines a state on $F(S)$.
\end{claim}
\begin{proof}Let $\widetilde S$ be a subspace of $\overline S$ generated by
$s$ and $u$. Let $v(\neq u)$ be a unit vector
in $S$ and put $P=[u]+[v]$. Then
$\widetilde S=P\oplus P^\bot$. Set
$\tilde w=v-\langle v,u\rangle u$ and let
$w=\frac{\tilde w}{\Vert\tilde w\Vert}$.
Similarly, let $\tilde z=u-\langle u,v\rangle v$
and put $z=\frac{\tilde z}{\Vert\tilde z\Vert}$.
Then $P=[v]\oplus[z]=[u]\oplus[w]$.  Define the
map
\begin{align*}
T:\widetilde S&\to \widetilde S\\
 P\oplus P^\bot&\to \widetilde S\\
 p+p'=\alpha v+\beta z+p'&\mapsto\alpha u+\beta w+p'.
\end{align*}
$T$ is a unitary operator on $\widetilde S$, that is $T$ is a
bijective linear operator satisfying
$$
\langle x,y\rangle=\langle Tx,Ty\rangle
$$
for all $x,y\in \widetilde S$.  \\

By the continuity of $T$ we can
extend it over $\overline S$. With a harmless abuse of notation let us
denote the extension again by $T$.
We now show that if $A$ is a subspace of
$S$, then  $\overline{TA}=T\overline A$.  Since
$T$ is continuous it follows immediately that
$T\overline A\subset\overline{TA}$.  Let
$x\in\overline{TA}$.  Then $x=\lim_{i\to\infty}x_i$
where $x_i\in TA$ for all $i\in\mathbb N$.  Let
$y_i\in A$ be such that $Ty_i=x_i$.  Then we have
\begin{align*}
\Vert x_i-x_j\Vert^2&=\langle Ty_i-Ty_j,Ty_i-Ty_j\rangle\\
                    &=\langle T(y_i-y_j),T(y_i-y_j)\rangle\\
                    &=\langle y_i-y_j,y_i-y_j\rangle\\
                    &=\Vert y_i-y_j\Vert^2.
\end{align*}
This implies that $\{y_i\}$ is Cauchy and therefore
it converges to some $y\in\overline A$.  That $Ty=x$
follows again by the continuity of $T$.\\

We now show that for any $A\in F(S)$, we have
$$
\Vert P_{\overline A}v\Vert^2=\Vert P_{\overline{TA}}u\Vert^2.
$$
Let $\{a_i\}\subset A$ be an ONB of $\overline A$.
Then $\{Ta_i\}$ is an ONB of $T\overline A$
($=\overline{TA}$) in $TA$.  We then have
\begin{align*}
a_i&=\alpha_iv+\beta_iz+p_i'\qquad\text{and therefore}\\
Ta_i&=\alpha_iu+\beta_iw+p_i'.
\end{align*}
This implies that
\begin{align*}
\Vert P_{\overline A}v\Vert^2&=\sum{|\langle a_i,v\rangle|^2}\\
                             &=\sum{|\alpha_i|^2}\\
                             &=\sum{|\langle Ta_i,u\rangle|^2}\\
                             &=\Vert P_{\overline{TA}}u\Vert^2.
\end{align*}
Thus, for any $A\in F(S)$, $s_v(A)=
\Vert P_{\overline A}v\Vert^2=\Vert P_{\overline{TA}}u\Vert^2=
s_u(TA)$, and therefore $s_v$ does indeed define a state on $F(S)$.
\end{proof}

\begin{claim}Suppose that, for each $u\in S$, $s_u$ defines a state on
$F(S)$. Then for every unit vector $v\in\overline S$,
$s_v$ defines a state on $F(S)$.
\end{claim}
\begin{proof}
Let $v\in\overline S\setminus S$.  There exists a
sequence $\{v_i\}\subset S$ such that $v=\lim_{i\to\infty}v_i$.
For any $A\in F(S)$,
\begin{align*}
P_{\overline A}v&=P_{\overline A}\lim_{i\to\infty}v_i\\
                &=\lim_{i\to\infty}P_{\overline A}v_i
\end{align*}
and therefore
$$
s_v(A)=\lim_{i\to\infty}s_{v_i}(A).
$$
It is then not difficult to check that $s_v$ defines
a state on $F(S)$ (pointwise limits of finitely additive states are finitely
additive states).
\end{proof}

\begin{claim}Let for any $v \in \overline S$ $s_v$ defines a state on
$F(S)$. Let $M$ be a closed subspace of $S$.
Then
$$
M\in F(S)\text{ if, and only if, }\overline{M^{\bot_S}}=\overline{M}^{\bot_{\overline S}}.
$$
\end{claim}
\begin{proof}
Let $M\in F(S)$.  We need to show that
$\overline{M^{\bot_S}}=\overline{M}^{\bot_{\overline S}}$.
It is sufficient to prove that $\overline{M^{\bot_S}}\supset\overline {M}^{\bot_{\overline S}}$.
Let $\{n_i:i\in I_{M'}\}$ be an orthonormal basis (ONB) in $M^{\bot_S}$ of $\overline{M^{\bot_S}}$,
and let $\tilde x\in \overline M^{\bot_{\overline S}}$ ($\tilde x\neq 0$) be arbitrary.
Put $x=\frac{\tilde x}{\Vert\tilde x\Vert}$.  Consider the state $s_x$ on $F(S)$.
\begin{align*}
1=s_x(S)&=s_x(M\vee M^{\bot_S})\\
        &=s_x(M)+s_x(M^{\bot_S})\\
        &=s_x(M^{\bot_S})\qquad \text{since }x\bot M\\
        &=\sum_{i\in I_{M'}}|\langle x,n_i\rangle|^2.
\end{align*}
This implies that for all $\tilde x\in \overline M^{\bot_{\overline S}}$,
$$
\Vert\tilde x\Vert^2=\sum_{i\in I_{M'}}|\langle\tilde x,n_i\rangle|^2.
$$
Therefore it follows, by Parseval's identity, that $\{n_i:i\in I_{M'}\}$
is an ONB of $\overline M^{\bot_{\overline S}}$ and hence
$\overline{M^{\bot_S}}=\overline M^{\bot_{\overline S}}$.\\
Now we prove the converse.  Suppose that
$\overline{M^{\bot_S}}=\overline{M}^{\bot_{\overline S}}$.  To reach a
contradiction, assume that
$M\notin F(S)$.  There exists $v\in M^{\bot_S\bot_S}\setminus M$
such that $v\bot M^{\bot_S}$ and $v\notin M$.  This implies that
$v\bot\overline{M^{\bot_S}}$ and hence
$v\in \overline{M}^{\bot_{\overline S}\bot_{\overline S}}=\overline M$.
But this would imply that $v\in\overline M\cap S=M$, since
$M$ is closed in $S$.  This is the required contradiction.
\end{proof}

\begin{claim}Suppose that for every $u \in S$ the mapping
  $$M\mapsto\langle P_{\overline M}u,u\rangle$$
defines a state on $F(S)$.
Let $M\in F(S)$ and let $\{x_i\}$ be any maximal orthonormal system (MONS) in $M$.  Then
$M=\{x_i\}^{\bot_S\bot_S}$.
\end{claim}
\begin{proof}
Let $\{m_i\}\subset M$ be an ONB of $\overline M$ and $\{n_i\}\subset M^{\bot_S}$
be an ONB of $\overline{M^{\bot_S}}=\overline{M}^{\bot_{\overline S}}$.  Then
$\{x_i\}\cup\{n_i\}$ is a MONS of $S$.  This implies that
\begin{align*}
&M\vee M^{\bot_S}\\
=&\vee[m_i]\bigvee\vee[n_i]\\
=&S\\
                                          =&\vee[x_i]\bigvee\vee[n_i]\\
                                          =&\{x_i\}^{\bot_S\bot_S}\vee M^{\bot_S}.
\end{align*}
Certainly, we have $\{x_i\}^{\bot_S\bot_S}\subset M$. Take any
unit vector $y\in \overline{M}$ and
consider the state $s_y$.  We have
\begin{align*}
1=s_y(S)&=s_y(\{x_i\}^{\bot_S\bot_S}\vee M^{\bot_S})\\
        &=s_y(\{x_i\}^{\bot_S\bot_S})\\
        &=\Vert P_{\overline{\{x_i\}^{\bot_S\bot_S}}}y\Vert.
\end{align*}
This implies that $y\in\overline{\{x_i\}^{\bot_S\bot_S}}$ and therefore
$$
\overline M=\overline{\{x_i\}^{\bot_S\bot_S}}
$$
which yields
$$
M=\{x_i\}^{\bot_S\bot_S}.
$$
\end{proof}

\begin{claim}$F(S)$ is orthomodular.
\end{claim}
\begin{proof}
Let $A\subset B$ be in $F(S)$.  Let $\{a_i\}\subset A$ be an ONB of
$\overline A$. Extend $\{a_i\}$ to a MONS $\{a_i\}\cup\{b_i\}$ of $B$.
%Let $\{c_i\}\subset B^{\bot_S}$ be an ONB of $\overline{B^{\bot_S}}$.
%Then $\{c_i\}^{\bot_S}=(\{a_i\}\cup\{b_i\})^{\bot_S}$.
It is not difficult to see that $\{b_i\}$ is a MONS in $A^{\bot_S}\cap B$
and that therefore
\begin{align*}
A\vee(A^{\bot_S}\wedge
B)&=\{a_i\}^{\bot_S\bot_S}\vee\{b_i\}^{\bot_S\bot_S}\\
&=\vee[a_i]\bigvee\vee[bi]\\
&=B
\end{align*}
\end{proof}
This completes the proof of Theorem \ref{t:1}.
\end{proof}

\providecommand{\bysame}{\leavevmode\hbox to3em{\hrulefill}\thinspace}

\end{document}